\documentclass[]{amsart}

\usepackage{latexsym}
\usepackage{amssymb,amsmath,amsopn}
\usepackage[dvips]{graphicx}   
\usepackage{color,epsfig}      
\usepackage{url}

\newtheorem{theorem}{Theorem}[section]
\newtheorem{lemma}{Lemma}[section]
\newtheorem{statement}{Statement}[section]
\newtheorem{proposition}{Proposition}[section]
\newtheorem{defi}{Definition}[section]
\newtheorem{problem}{Problem}[section]

\theoremstyle{definition}
\newtheorem{corollary}{Corollary}[section]
\newtheorem{conjecture}{Conjecture}[section]
\newtheorem{remark}[]{Remark}[section]

\DeclareMathOperator{\conv}{conv}

\DeclareMathOperator{\vol}{vol}
\DeclareMathOperator{\area}{area}
\DeclareMathOperator{\inter}{int}
\DeclareMathOperator{\relint}{relint}

\DeclareMathOperator{\bd}{bd}

\newcommand{\K}{\mathcal{K}}
\renewcommand{\S}{\mathcal{S}}

\newcommand{\B}{\mathbf B}
\newcommand{\Sph}{\mathbb{S}}
\renewcommand{\P}{\mathcal{P}}

\DeclareMathOperator{\Sym}{Sym}

\DeclareMathOperator{\aff}{aff}

\renewcommand{\S}{\mathcal{S}}
\renewcommand{\Re}{\mathrm{ Re}}
\newcommand{\R}{\mathbb{R}}
\newcommand{\Z}{\mathbb{Z}}

\renewcommand{\P}{\mathcal{P}}
\newcommand{\F}{\mathcal{F}}
\newcommand{\C}{\mathcal{C}}

\title{Volume of convex hull of two bodies and related problems }
\author{\'Akos G.Horv\'ath}
\address{\'Akos G.Horv\'ath, Dept. of Geometry, Budapest University of Technology,
Egry J\'ozsef u. 1., Budapest, Hungary, 1111}
\email{ghorvath@math.bme.hu}


\begin{document}

\begin{abstract}
In this paper we deal with problems concerning the volume of the convex hull of two "connecting" bodies. After a historical background we collect some results, methods and  open problems, respectively.	
\end{abstract}

\maketitle

\section{Introduction}

To find the convex polyhedra in Euclidean $3$-space $\mathbb{R}^3$, with a given number of faces and with minimal isoperimetric quotient, is a centuries old question
of geometry: research in this direction perhaps started with the work of Lhuilier in the 18th century. A famous result of Lindel\"of \cite{L70}, published in the 19th century, yields a necessary condition for such a polyhedron: it states that any optimal polyhedron is circumscribed about a Euclidean ball, and this ball touches each face at its centroid. In particular, it follows from his result that, instead of fixing surface area while looking for minimal volume, we may fix the inradius of the polyhedron. Since the publication of this result, the same condition for polytopes in $n$-dimensional space $\mathbb{R}^n$ has been established (cf. \cite{G07}), and many variants of this problem have been investigated (cf., e.g. \cite{BB96}).
For references and open problems of this kind, the interested reader is referred to \cite{F93}, \cite{croft} or \cite{BMP05}. For polytopes with $(n+2)$ vertices this question was answered by Kind and Kleinschmidt \cite{KK76}.
The solution for polytopes with $n+3$ vertices was published in \cite{KW03}, which later turned out to be incomplete (cf. \cite{KW05}), and thus,
this case is still open.
We mention two problems in more detail:
\begin{itemize}
\item The dual of the original problem: to find, among $d$-polytopes with a given number of vertices and inscribed in the unit sphere, the ones with maximal volume, and

\item to find the extremity of the volume of the convex hull of two "connecting" bodies.
\end{itemize}

The first problem that to find the maximal volume polyhedra in $\mathbb{R}^3$ with a given number of vertices and inscribed in the unit sphere,
was first mentioned in \cite{ftl} in 1964.
A systematic investigation of this question was started with the paper \cite{bermanhanes} of Berman and Hanes in 1970, who found a necessary condition
for optimal polyhedra, and determined those with $n \leq 8$ vertices.
The same problem was examined in \cite{M02}, where the author presented the results of a computer-aided search for optimal polyhedra with $4 \leq n \leq 30$ vertices.
Nevertheless, according to our knowledge, this question, which is listed in both research problem books \cite{BMP05} and \cite{croft}, is still open for polyhedra with $n > 8$ vertices.

The second problem connected with the first one on the following way: If the given points form the respective vertex sets of two polyhedra (inscribed in the unit sphere) then the volume of the convex hull of these points is the volume of the  convex hull of two "connecting" bodies, too. It is interesting that the case of two regular simplices with common center gives another maximum as the global isodiametric problem on eight points inscribed in the unit sphere.

The examination of the volume of the convex hull of two congruent copies of a convex body in Euclidean $d$-space (for special subgroups) investigated systematically first by Rogers, Shepard and Machbeth in 1950s (see in \cite{ rogers-shephard 1}, \cite{ rogers-shephard2} and \cite{macbeath}). Fifty years later a problem similar to that of the simplices arose that lead to new investigations by new methods which obtained fresh results (see in \cite{gho-langi}, \cite{gho}, \cite{gho 2}).  In particular, a related conjecture of Rogers and Shephard has been proved in \cite{gho-langi}.

Finally we review some important consequences of the icosahedron inequality  of L. Fejes-T\'oth. In particular, it is needed for the proof of the statement that the maximal volume polyhedron spanned by the vertices of two regular simplices with common centroid is the cube. It is also used in the proof of that the maximal volume polyhedron with eight vertices and inscribed in the unit sphere is a triangular one distinct from the cube.

\section{Maximal volume polytopes inscribed in the unit sphere}

The aim of this section is to review the results on the first problem mentioned in the introduction.

Let for any $p,q \in \mathbb{R}^d$, $|p|$ and $[p,q]$ denote the standard Euclidean norm of $p$, and the closed segment with endpoints $p$ and $q$, respectively.
The origin of the standard coordinate system of $\mathbb{R}^d$ is denoted by $o$.
If $v_1,v_2,\ldots,v_d \in \mathbb{R}^d$, then the $d \times d$ determinant with columns $v_1,v_2,\ldots,v_d$,
in this order, is denoted by $|v_1,\ldots,v_d|$.
The unit ball of $\mathbb{R}^d$, with $o$ as its center, is denoted by $\mathrm{B}^d$, and we set $\Sph^{d-1}=\bd \B^d$.

Throughout this section, by a polytope we mean a convex polytope.
The vertex set of a polytope $P$ is denoted by $V(P)$.
We denote the family of $d$-dimensional polytopes, with $n$ vertices and inscribed in the unit sphere $\Sph^{d-1}$, by
$\P_d(n)$. The $d$-dimensional volume denotes by $\vol_d$, and set $v_d(n) = \max \{ \vol_d(P) : P \in \P_d(n) \}$.
Note that by compactness, $v_d(n)$ exists for any value of $d$ and $n$.

Let $P$ be a $d$-polytope inscribed in the unit sphere $\Sph^{d-1}$, and let $V(P) = \{ p_1,p_2,\ldots,p_n\}$.

Let $\C(P)$ be a simplicial complex with the property that $|\C(P)| = \bd P$, and that the vertices of $\C(P)$ are exactly the points of $V(P)$.
Observe that such a complex exist. Indeed, for any positive integer $k$, and for $i=1,2,\ldots,n$, consider a point $p^k_i$ such that
$|p^k_i - p_i| < \frac{1}{k}$, and the polytope $P_k=\conv \{ q^k_i : i=1,2,\ldots,n \}$ is simplicial.
Define $\C(P_k)$ as the family of the faces of $P_k$. We may choose a subsequence of the sequence $\{\C(P_k)\}$ with the property
that the facets of the complexes belong to vertices with the same indices.
Then the limit of this subsequence yields a complex with the required properties. Note that if $P$ is simplicial, then $\C(P)$ is the family of the faces of $P$.

Now we orient $\C(P)$ in such a way that for each $(d-1)$-simplex $(p_{i_1},p_{i_2},\ldots, p_{i_d})$ (where $i_1\leq i_2\leq \cdots \leq i_d$) in $C(P)$,
the determinant $|p_{i_1},\ldots,p_{i_d}|$ is positive; and call the $d$-simplex $\conv \{ o,p_{i_1},\ldots,p_{i_d}\}$ a \emph{facial simplex} of $P$.
We call the $(d-1)$-dimensional simplices of $\C(P)$ the \emph{facets} of $\C(P)$.

\subsection{$3$-dimensional results.}

The problem investigated in this section was raised by L. Fejes-T\'oth in \cite{ftl}. His famous inequality (called by icosahedron inequality) can be formulated as follows.

\begin{theorem}[\cite{ftl} on p.263]
If $V$ denotes the volume, $r$ the inradius and $R$ the circumradius of a convex polyhedron having $f$ faces, $v$ vertices and $e$ edges, then
\begin{equation}\label{ftineq}
\frac{e}{3}\sin\frac{\pi f}{e}\left(\tan^2\frac{\pi f}{2e}\tan^2\frac{\pi v}{2e}\right)r^3\leq V\leq \frac{2e}{3}\cos^2\frac{\pi f}{2e}\cot \frac{\pi v}{2e}\left(1-\cot^2\frac{\pi f}{2 e}\cot^2\frac{\pi v}{2 e}\right)R^3.
\end{equation}
Equality holds in both inequalities only for regular polyhedra.	
\end{theorem}

He noted that \emph{"a polyhedron with a given number of faces $f$ is always a limiting figure of a trihedral polyhedron with $f$ faces. Similarly, a polyhedron with a given number $v$ of vertices is always the limiting figure of a trigonal polyhedron with $v$ vertices. Hence introducing the notation
$$
\omega_n=\frac{n}{n-2}\frac{\pi}{6}
$$
we have the following inequalities
\begin{equation}\label{ftineqfaces}
(f-2)\sin 2\omega_f\left( 3\tan^2\omega_f-1\right) r^3\leq V\leq \frac{2\sqrt{3}}{9}\left( f-2\right) \cos^2\omega_f\left( 3-\cot^2\omega_f\right) R^3,
\end{equation}
\begin{equation}\label{ftineqvert}
\frac{\sqrt{3}}{2}(v-2)\left( 3\tan^2\omega_v-1\right) r^3\leq V\leq \frac{1}{6}\left( v-2\right) \cot\omega_v\left( 3-\cot^2\omega_v\right) R^3.
\end{equation}
Equality holds in the first two inequalities only for regular tetrahedron, hexahedron and dodecahedron (f=4, 6, 12) and in the last two inequalities only for the regular tetrahedron, octahedron and icosahedron (v=4, 6, 12)."}

The right hand side of inequality (1) immediately solves our first problem in the cases when the number of vertices is $v=4, 6, 12$; the maximal volume polyhedra with $4, 6$ and $12$ vertices inscribed in the unit sphere are the regular tetrahedron, octahedron and icosahedron, respectively.

The second milestone in the investigation of this problem is the paper of Berman and Hanes (\cite{bermanhanes}) written in 1970. They solved the problem for  $v=5, 7, 8$ vertices, respectively. Their methods are based on a combinatorial classification of the possible spherical tilings due to Bowen and Fisk (\cite{bowenfisk}) and a geometric result which gives a condition for the local optimal positions. They characterized these positions by a property called Property Z. We now give the definitions with respect to the $d$-dimensional space.

\begin{defi}
	Let $P \in \P_d(n)$ be a $d$-polytope with $V(P) = \{p_1,p_2,\ldots ,p_n\}$.
	If for each $i$, there is an open set $U_i\subset \Sph^{d-1}$ such that $p_i\in U_i$, and for any $q \in U_i$, we have
	\[
	\vol_d\left( \conv \left( \left( V(P) \setminus \{p_i \} \right) \cup \{ q \} \right)\right) \leq \vol_d\left(P\right) ,
	\]
	then we say that $P$ satisfies \emph{Property Z}.
\end{defi}

Returning to the three-dimensional case if $p_i$ and $p_j$ are vertices of $P$, denote the line segment whose endpoints are $p_i$ and $p_j$ by $s_{ij}$ and its length by $|s_{ij}|$. Also, let $n_{ij}=1/6\left( p_i\times p_j\right) $ where $\times$ denotes the vector product in $E^3$.
\begin{lemma}[Lemma 1 in \cite{bermanhanes}]\label{propz3}
Let $P$ with vertices $p_1,\ldots,p_n$ have property $Z$. Let $C(P)$ be any oriented complex associated with $P$ such that $\vol(C(P))\geq 0$. Suppose $s_{12},\ldots ,s_{1r}$ are all the edges of $C(P)$ incident with $p_1$ and that $p_2, p_3, p_1$; $p_3, p_4, p_1$; ... ;$p_r, p_2, p_1$ are orders for faces consistent with the orientation of $C(P)$.
\begin{itemize}
	\item[i,] Then $p_1=m/|m|$ where $m=n_{23}+n_{34}+\cdots +n_{r2}$.
	\item[ii,] Furthermore, each face of $P$ is triangular.
\end{itemize} 	
\end{lemma}
Let the \emph{valence} of a vertex of $C(P)$ be the number of edges of $C(P)$ incident with that vertex. By Euler's formula the average of the valences is $6-12/n$. If $n$ is such that $6-12/n$ is an integer then $C(P)$ is \emph{medial} if the valence of each vertex is $6-12/n$. If $6-12/n$ is not an integer then $C(P)$ is \emph{medial} provided the valence of each vertex is either $m$ or $m+1$ where $m<6-12/n<m+1$. $P$ is said to be \emph{medial} provided all faces of $P$ are triangular and $C(P)$ is medial. Goldberg in \cite{goldberg} made a conjecture whose dual was formulated by Grace in \cite{grace}: The polyhedron with $n$ vertices in the unit sphere whose volume is a maximum is a medial polyhedron provided a medial polyhedron exists for that $n$.
Connecting to this conjecture  Berman and Hanes proved that if $n=4,5,6,7,8$ then the polyhedra with maximal volume inscribed in the unit sphere are medial polyhedra with Property Z, respectively.

Note that in the proofs of the above results (on $n\geq 5$) is an important step to show that the valences of the vertices of a polyhedron with maximal volume are at least $4$. This follows from inequality (2).

\begin{figure}[h]
\includegraphics[scale=1]{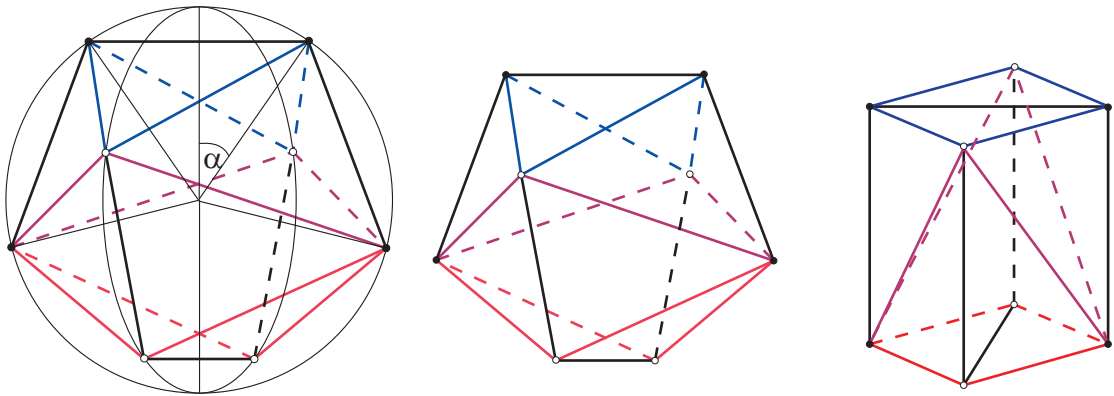}
\caption{The medial complex with $8$ vertices and its two polyhedra, the maximal volume polyhedron and the cube.}
\end{figure}

The maximal volume polyhedron for $n=4$ is the regular simplex. For $n=5,6,7$ they are the so-called double $n$-pyramids, with $n=5,6,7$, respectively. (By a \emph{double $n$-pyramid} (for $n\geq 5$), is meant a complex of $n$ vertices with two vertices of valence $n-2$ each of which is connected by an edge to each of the remaining $n-2$ vertices, all of which have valence $4$. The $2(n-2)$ faces of a double $n$-pyramid are all triangular. A polyhedron $P$ is a \emph{double $n$-pyramid} provided each of its faces is triangular and some $C(P)$ is a double $n$-pyramid.) An interesting observation (see Lemma 2 in \cite{bermanhanes}) is that if $P$ is a double $n$-pyramid with property Z then $P$ is unique up to congruence and its volume is $[(n-2)/3]\sin 2\pi/(n-2)$.

For $n=8$ there exists only two non-isomorphic complexes which have no vertices of valence $3$ (see in \cite{bowenfisk}). One of them the double $8$-pyramid and the other one has four valence 4 vertices and four valence 5 vertices, and therefore it is the medial complex (see on Fig.1). It has been shown that if this latter has Property Z then $P$ is uniquely determined up to congruence and its volume is
$
\sqrt{\left[\frac{475+29\sqrt{145}}{250}\right]}
$
giving the maximal volume polyhedron with eight vertices.

As concluding remarks Berman and Hanes raised the following questions:
\begin{problem}
	 For which types of polyhedra does Property Z determine a unique polyhedron. More generally, for each isomorphism class of polyhedra is there one and only one polyhedron (up to congruence) which gives a relative maximum for the volume?
\end{problem}

\begin{problem}
	 For $n=4,\ldots, 7$ the duals of the polyhedra of maximum volume are just those polyhedra with $n$ faces circumscribed about the unit sphere of minimum volume. For $n=8$ the dual of the maximal volume polyhedron (described above) is the best known solution to the isoperimetric problem for polyhedra with $8$ faces. Is this true in general?
\end{problem}

Recently there is no answer for these questions.

We have to mention a theorem of A. Florian (\cite{florian})which immediately implies the inequalities in (3). Let $P$ be a convex polyhedron with $v$ vertices and volume $\vol(P)$. We consider an orthoscheme $T=OABC$ (where $OA$ orthogonal to the plane $ABC$, and $AB$ orthogonal to $BC$) with the properties:
\begin{itemize}
	\item[(i)] the radial projection of $ABC$ onto the unit sphere with centre $O$ is the spherical triangle $T'=A'B'C'$ given by
	$$
	C'A'B'\sphericalangle=\frac{\pi}{3}, \quad A'B'C'\sphericalangle =\frac{\pi}{2}, \quad \area(T')= \frac{4\pi}{12(v-2)}
	$$
	\item[(ii)]
	$
	\vol(T)=\frac{1}{12(v-2)}\vol(P).
	$
\end{itemize}
Then we have:
\begin{theorem}[\cite{florian}]
	Let $K(\rho)$ be a ball with centre $O$ and radius $\rho$. Let $P$ be a convex polyhedron with $v$ vertices and volume $\vol(P)$, and let the tetrahedron $T$ be defined as above. Then
	\begin{equation}
	\vol(P\cap K(\rho))\leq 12(v-2)\vol(T\cap K(\rho))
	\end{equation}
	with equality if $v=4, 6$ or $12$ and $P$ is a regular tetrahedron, octahedron or icosahedron with centre $O$. When $|OA|\leq \rho \leq |OC|$, these are the only cases of equality.
\end{theorem}

We recall the paper of Mutoh \cite{M02} who presented the results of a computer-aided search for optimal polyhedra with $4 \leq n \leq 30$ vertices. The solutions of the computation probably solved the mentioned cases, respectively, however there is no information in the paper either on the source code of the program or the algorithm  which based the computation.
Table \ref{mutohtable}  describes some of those polyhedra which suggested by the author as the maximal volume one inscribed in the unit sphere. We refer here only a part of the complete table of Mutoh, for more information see the original paper \cite{M02}.

\begin{table}[ht]
  \centering
\begin{tabular}{|c|c|c|c|}
\hline
Number of vertices & Maximal volume & Number of Facets & Valences of vertices \\
\hline
4 & 0.51320010 & 4 & $3 \times 4$ \\
\hline
5 & 0.86602375 & 6 & $3\times 2$, $4\times 3$ \\
\hline
6 & 1.33333036 & 8 & $4 \times 6$ \\
\hline
7 & 1.58508910 & 10 & $4\times 5$, $5\times 2$ \\
\hline
8 & 1.81571182 & 12 & $4\times 4$, $5\times 4$ \\
\hline
9 & 2.04374046 & 14 & $4\times 3$, $5\times 6$ \\
\hline
10 & 2.21872888 & 16 & $4\times 2$, $5\times 8$ \\
\hline
11 & 2.35462915 & 18 & $4\times 2$, $5\times 8$, $6\times 1$ \\
\hline
12 & 2.53614471 & 20 & $5\times 12$ \\
\hline
\vdots & \vdots & \vdots & \vdots \\
\hline
30 & 3.45322727 & 56 & $5\times 12$, $6\times 18$ \\
\hline
\end{tabular}
\vspace{0.2cm}
  \caption{Computer search results of polyhedra of maximal volume inscribed in the unit sphere.}\label{mutohtable}
\end{table}
Mutoh notes that it seems to be that the conjecture of Grace on medial polyhedron is false because the polyhedra are just found by computer optimal in the cases $n=11$ and $n=13$ are not medial ones, respectively. Mutoh also listed the polyhedra circumscribed to the unit sphere with minimal volume and examined the dual conjecture of Goldberg (see also Problem 2.2).
He said:
\emph{"Goldberg conjectured that the polyhedron of maximal volume inscribed to
the unit sphere and the polyhedron of minimal volume circumscribed about the
unit sphere are dual. A comparison of Table 1 and 3 shows that the number
of vertices and the number of faces of the two class of polyhedra correspond
with each other. The degrees of vertices of the polyhedra of maximal volume
inscribed in the unit sphere correspond to the numbers of vertices of faces of
the polyhedra of minimal volume circumscribed about the unit sphere. Indeed,
the volume of polyhedra whose vertices are the contact points of the unit sphere
and the polyhedra circumscribed about the unit sphere differs only by 0.07299\%
from the volume of the polyhedra inscribed in the unit sphere."}

We turn to a recent result that generalizes the triangle case of the inequality (\ref{ftineq}) of L. Fejes-T\'oth. If $A,B,C$ are three points on the unit sphere we can consider two triangles, one of the corresponding spherical triangle and the second one the rectilineal triangle with these vertices, respectively. Both of them are denoted by $ABC$. The angles of the rectilineal triangle are the half of the angles between those radius of the circumscribed circle which connect the center $K$ of the rectilineal triangle $ABC$ to the vertices $A,B,C$. Since $K$ is also the foot of the altitude of the tetrahedron with base $ABC$ and apex $O$, hence the angles $\alpha_A$, $\alpha_B$ and $\alpha_C$ of the rectilineal triangle $ABC$, play an important role in our investigations, we refer to them as the \emph{central angles} of the spherical edges $BC$, $AC$ and $AB$, respectively. We call again the tetrahedron $ABCO$ the \emph{facial tetrahedron} with base $ABC$ and apex $O$.

\begin{lemma}(See in \cite{gho 3})\label{isosceles}
	Let $ABC$ be a triangle inscribed in the unit sphere. Then there is an isosceles triangle $A'B'C'$ inscribed in the unit sphere with the following properties:
	\begin{itemize}
		\item the greatest central angles and also the spherical areas of the two triangles are equal to each other, respectively;
		\item the volume of the facial tetrahedron with base $A'B'C'$  is greater than or equal to the volume of the facial tetrahedron with base $ABC$.
	\end{itemize}
\end{lemma}
From Lemma \ref{isosceles} it can be proved upper bound functions for the volume of the facial tetrahedron.
\begin{proposition}
	Let the spherical area of the spherical triangle $ABC$ be $\tau$. Let $\alpha_C$ be the greatest central angle of $ABC$ corresponding to $AB$. Then the volume $V$ of the facial tetrahedron $ABCO$ holds the inequality
	\begin{equation}
	V\leq \frac{1}{3}\tan\frac{\tau}{2}\left(2-\frac{|AB|^2}{4}\left(1+\frac{1}{\left(1+\cos\alpha_C\right)}\right)\right).
	\end{equation}
	In terms of $\tau $ and $c:=AB$ we have
	\begin{equation}
	V\leq v(\tau,c):=\frac{1}{6}\sin c\frac{\cos\frac{\tau-c}{2}-\cos\frac{\tau}{2}\cos\frac{c}{2}}{1-\cos\frac{c}{2}\cos\frac{\tau}{2}}.
	\end{equation}
	Equality holds if and only if $|AC|=|CB|$.
\end{proposition}

Observe that the function $v(\tau,c)$ is concave in the parameter domain $\mathcal{D}:=\{0<\tau<\pi/2, \tau\leq c<\min\{f(\tau),2\sin^{-1}\sqrt{2/3}\}\}$ with certain concave (in $\tau$) function $f(\tau)$ defined by the zeros of the Hessian; and non-concave in the domain $\mathcal{D'}=\{0< \tau \leq \omega, f(\tau)\leq c\leq 2\sin^{-1}\sqrt{2/3}\}=\{0< \tau\leq c\leq \pi/2\}\setminus D$, where $f(\omega)=2\sin^{-1}\sqrt{2/3}$.

Assume now that the triangular star-shaped polyhedron $P$ with $f$ face inscribed in the unit sphere.  Let $c_1,\ldots,c_f$ be the arc-lengths of the edges of the faces $F_1,\ldots ,F_f$ corresponding to their maximal central angles, respectively. Denote by $\tau_i$ the spherical area of the spherical triangle corresponding to the face $F_i$ for all $i$. We note that for a spherical triangle which edges $a,b,c$ hold the inequalities $0<a\leq b\leq c<\pi/2$, also holds the inequality $\tau \leq c$. In fact, for fixed $\tau$ the least value of the maximal edge length attend at the case of regular triangle. If $c<\pi/2$ then we have
$$
\tan\frac{\tau}{4} =\left(\tan\frac{c}{4}\sqrt{\tan\frac{3c}{4}\tan\frac{c}{4}}\right)= \left(\tan\frac{c}{4}\sqrt{1-\frac{\tan\frac{3c}{4}+\tan\frac{c}{4}}{\tan c }}\right)<\tan\frac{c}{4},
$$
and if $c=\pi/2$ then $\tau=8\pi/4=\pi/2$ proving our observation.

The following theorem gives an upper bound on the volume of the star-shaped polyhedron corresponding to the given spherical tiling in question.

\begin{theorem}(See in \cite{gho 3})
Assume that $0<\tau_i<\pi/2$ holds for all $i$. For $i=1,\ldots, f'$ we require the inequalities $0<\tau_i\leq c_i\leq \min\{f(\tau_i),2\sin^{-1}\sqrt{2/3}\}$ and for all $j$ with $j\geq f'$ the inequalities $0< f(\tau_j)\leq c_j \leq 2\sin^{-1}\sqrt{2/3}$, respectively. Let denote $c':=\frac{1}{f'}\sum\limits_{i=1}^{f'}c_i$, $c^\star:=\frac{1}{f-f'}\sum\limits_{i=f'+1}^{f}f(\tau_i)$ and $\tau':=\sum\limits_{i=f'+1}^{f}\tau_i$, respectively. Then we have
	\begin{equation}
	v(P)\leq \frac{f}{6}\sin \left(\frac{f'c'+(f-f')c^\star}{f}\right)\frac{\cos \left(\frac{4\pi-f'c'-(f-f')c^\star}{2f}\right)-\cos\frac{2\pi}{f}\cos\left(\frac{f'c'+(f-f')c^\star}{2f}\right)} {1-\cos\frac{4\pi}{2f}\cos\left(\frac{f'c'+(f-f')c^\star}{2f}\right)} .
	\end{equation}
\end{theorem}

\subsection{The cases of higher dimensions}

As we saw in the previous subsection, even the $3$-dimensional case is completely proved only when the number of vertices less or equal to eight.  This shows that in higher dimensions we cannot expect such complete results as was published by L. Fejes-T\'oth, A. Fl\'orian or Berman and Hanes in the second half on the last century, respectively. As Fl\'orian said in \cite{florian}: "Several extremum properties of the regular triangle and the regular tetrahedron may be generalized to regular simplices in all dimensions.... Little is known in this respect about the general cross polytope, the hypercube and the nontrivial regular convex polytopes in $4$-space". Some extremum properties of these polytopes were established by comparing them with the topologically isomorphic convex polytopes....But no methods are available for proving inequalities analogous to (\ref{ftineqvert})." We now extract the method of Berman and Hanes to higher dimensions and using a combinatorial concept, the idea of Gale's transform
solve some cases of few vertices. In this subsection we collect the results of the paper \cite{gho-langi 2}.

The first step is the generalization of Lemma \ref{propz3} for arbitrary dimensions.
\begin{lemma}\label{lem:propz}
Consider a polytope $P \in \P_d(n)$ satisfying Property Z.
For any $p \in V(P)$, let $\F_p$ denote the family of the facets of $\C(P)$ containing $p$.
For any $F \in \F_p$, set
$$
A(F,p) = \vol_{d-1} \left( \conv \left( \left( V(F) \cup \{ o \} \right) \setminus \{ p \} \right) \right),
$$ and let $m(F,p)$ be the unit normal vector of the hyperplane, spanned by $\left( V(F) \cup \{ o \} \right) \setminus \{p\}$, pointing in the direction of the half space containing $p$.
		\begin{enumerate}
		\item[(\ref{lem:propz}.1)] Then we have
		$
		p = m/|m|, \mbox{ where } m=\sum_{F \in \F_p} A(F,p) m(F,p).
		$
		\item[(\ref{lem:propz}.2)] Furthermore $P$ is simplicial.
	\end{enumerate}
\end{lemma}

\begin{remark}\label{rem:hyperplane}
	Assume that $P \in \P_d(n)$ satisfies Property Z, and for some $p \in V(P)$, all the vertices of $P$ adjacent to $p$ are contained in a hyperplane $H$.
	Then the supporting hyperplane of $\Sph^{d-1}$ at $p$ is parallel to $H$, or in other words, $p$ is a normal vector to $H$.
	Thus, in this case all the edges of $P$, starting at $p$, are of equal length.
\end{remark}

\begin{lemma}\label{lem:note2}
	Let $P \in \P_d(n)$ satisfy Property Z, and let $p \in V(P)$. Let $q_1, q_2 \in V(P)$ be adjacent to $p$.
	Assume that any facet of $P$ containing $p$ contains at least one of $q_1$ and $q_2$, and for any $S \subset V(P)$ of cardinality $d-2$, $\conv (S \cup \{p,q_1\})$ is a facet of $P$ not containing $q_2$ if, and only if $\conv (S \cup \{p,q_2\})$ is a facet of $P$ not containing $q_1$. Then $|q_1-p| = |q_2-p|$.
\end{lemma}

Corollary~\ref{cor:simplex} is a straightforward consequence of Lemma~\ref{lem:note2} or, equivalently, Remark~\ref{rem:hyperplane}.

\begin{corollary}\label{cor:simplex}
	If $P \in \P_d(d+1)$ and $\vol_d(P) = v_d(d+1)$, then $P$ is a regular simplex inscribed in $\Sph^{d-1}$.
\end{corollary}

We note that this statement can be considered as a folklore. The analogous statement in $d$-dimensional spherical geometry (for simplices inscribed in a sphere of $\S^d$ with radius less than $\pi/2$) was proved by K. B\"or\"oczky in \cite{boroczky}. The method of B\"or\"oczky is based on the fact that Steiner's symmetrization is a volume-increasing transformation of the spherical space and so it can not be transformed immediately to the hyperbolic case. In hyperbolic spaces the investigations concentrated only to the simplices with ideal vertices. In dimension two every two triangles with ideal vertices are congruent to each other implying that they have the same area which value is maximal one  among the triangles. On the other hand it was proved by Milnor (see in \cite{milnor-thurston} or in \cite{milnor}) that in hyperbolic 3-space, a simplex is of maximal volume if and only if it is ideal and regular. The same $d$-dimensional statement has been proved by U. Haagerup and H. J. Munkholm in \cite{haagerupmunkholm}. This motivates the following:

\begin{problem}
Prove or disprove that in hyperbolic $d$-space a simplex is of maximal volume inscribed in
the unit sphere if and only if it is a regular one.
\end{problem}

Before the next corollary recall that if $K$ is a $(d-1)$-polytope in $\R^d$, and $[p_1,p_2]$ is a segment intersecting the relative interior of $K$
at a singleton different from $p_1$ and $p_2$, then $\conv (K \cup [p_1,p_2])$ is a \emph{$d$-bipyramid} with base $K$ and apexes $p_1,p_2$ (cf. \cite{grunbaum}). In the literature the terminology "bipyramid" is more prevalent as of the nomenclature "double-pyramid" of Berman and Hanes. In the rest of this paper we use bipyramid.

\begin{corollary}\label{cor:bipyramid}
	Let $P \in \P_d(n)$ be combinatorially equivalent to a $d$-bipyramid. Assume that $P$ satisfies Property Z.
	Then $P$ is a $d$-bipyramid, its apexes $p_1,p_2$ are antipodal points, its base $K$ and $[p_1,p_2]$ lie in orthogonal linear subspaces of $\R^d$, and $K$ satisfies Property Z in the hyperplane $\aff K$.
\end{corollary}

Corollary~\ref{cor:bipyramid} implies the following one:

\begin{corollary}\label{cor:crosspolytope}
	If $P \in \P_{d}(2d)$ has maximal volume in the combinatorial class of cross-polytopes inscribed in $\Sph^{d-1}$, then it is a regular cross-polytope.
\end{corollary}

The first non-trivial case is when the number of points is equal to $n=d+2$. It has been proved:

\begin{theorem}[\cite{gho-langi 2}]\label{thm:dplus2}
	Let $P \in \P_d(d+2)$ have maximal volume over $\P_d(d+2)$.
	Then $P=\conv (P_1 \cup P_2)$, where $P_1$ and $P_2$ are regular simplices of dimensions $\lfloor \frac{d}{2} \rfloor$ and $\lceil \frac{d}{2} \rceil$, respectively, inscribed in $\Sph^{d-1}$, and contained in orthogonal linear subspaces of $\R^d$.
	Furthermore,
	\[
	v_d(d+2) = \frac{1}{d!} \cdot \frac{\left( \lfloor d/2 \rfloor + 1 \right)^{\frac{\lfloor d/2 \rfloor + 1}{2}} \cdot \left( \lceil d/2 \rceil + 1 \right)^{\frac{\lceil d/2 \rceil + 1}{2}}}{\lfloor d/2 \rfloor^{\frac{\lfloor d/2 \rfloor}{2}} \cdot \lceil d/2 \rceil^{\frac{\lceil d/2 \rceil}{2} }}
	\]
\end{theorem}

In the proof of the results on $d$-polytopes with $d+2$ or $d+3$ vertices, we use extensively the properties of the so-called Gale transform of a polytope
(cf. \cite{grunbaum}, \cite{ziegler}). Since the application of this combinatorial theory leads to a new method in the investigation of our problem we review it.

Consider a $d$-polytope $P$ with vertex set $V(P)=\{p_i:i=1,2,\ldots,n \}$. Regarding $\R^d$ as the hyperplane $\{ x_{d+1} = 1\}$ of $\R^{d+1}$, we can represent $V(P)$ as a $(d+1) \times n$ matrix $M$, in which each column lists the coordinates of a corresponding vertex in the standard basis of $\R^{d+1}$. Clearly, this matrix has rank $d+1$, and thus, it defines a linear mapping $L : \R^n \rightarrow \R^{d+1}$, with $\dim \ker L= n-d-1$.
Consider a basis $\{ w_1,w_2,\ldots, w_{n-d-1} \}$ of $\ker L$, and let $\bar{L}:\R^{n-d-1} \rightarrow \R^n$ be the linear map mapping the $i$th vector of the standard basis of $\R^{n-d-1}$ into $w_i$.
Then the matrix $\bar{M}$ of $\bar{L}$ is an $n\times (n-d-1)$ matrix of (maximal) rank $n-d-1$, satisfying the equation $M \bar{M} = O$, where $O$ is the matrix with all entries equal to zero. Note that the rows of $\bar{M}$ can be represented as points of $\R^{n-d-1}$.
For any vertex $p_i \in V(P)$, we call the $i$th row of $\bar{M}$ the \emph{Gale transform of $p_i$}, and denote it by $\bar{p}_i$.
Furthermore, the $n$-element multiset $\{ \bar{p_i} : i=1,2,\ldots, n \} \subset \R^{n-d-1}$ is called the \emph{Gale transform of $P$}, and is denoted by $\bar{P}$. If $\conv S$ is a face of $P$ for some $S \subset V(P)$, then the (multi)set of the Gale transform of the points of $S$ is called  a face of $\bar{P}$.
If $\bar{S}$ is a face of $\bar{P}$, then $\bar{P} \setminus \bar{S}$ is called a \emph{coface} of $\bar{P}$.

Let $V = \{ q_i: i=1,2,\ldots,n\} \subset \R^{n-d-1}$ be a (multi)set. We say that $V$ is a \emph{Gale diagram} of $P$, if for some Gale transform $P'$
the conditions $o \in \relint \conv \{ q_j : j \in I \}$ and $o \in \relint \conv \{ \bar{p}_j : j \in I \}$ are satisfied for the same subsets of $\{1,2,\ldots,n\}$. If $V \subset \Sph^{n-d-2}$, then $V$ is a \emph{normalized Gale diagram} (cf. \cite{lee}). A \emph{standard Gale diagram} is a normalized Gale diagram in which the consecutive diameters are equidistant. A \emph{contracted Gale diagram} is a standard Gale diagram which has the least possible number of diameters among all isomorphic diagrams. We note that each $d$-polytope with at most $d+3$ vertices may be represented by a contracted Gale diagram (cf. \cite{grunbaum} or \cite{ziegler}).
An important tool of the proofs the following theorem from \cite{grunbaum} or also from \cite{ziegler}.

\begin{theorem}[\cite{grunbaum},\cite{ziegler}]\label{thm:Gale}
	\begin{itemize}
		\item[(i)] A multiset $\bar{P}$ of $n$ points in $\R^{n-d-1}$ is a Gale diagram of a $d$-polytope
		$P$ with $n$ vertices if and only if every open half-space in $\R^{n-d-1}$ bounded by a hyperplane
		through $o$ contains at least two points of $\bar{V}$ (or, alternatively, all the points of $\bar{P}$
		coincide with $o$ and then $n=d+1$ and $P$ is a  $d$-simplex).
		\item[(ii)] If $F$ is a facet of $P$, and $Z$ is the corresponding coface, then in any Gale
		diagram $\bar{V}$ of $P$, $\bar{Z}$ is the set of vertices of a (non-degenerate) set with $o$ in its
		relative interior.
		\item[(iii)] A polytope $P$ is simplicial if and only if, for every hyperplane $H$ containing
		$o  \in \R^{n-d-1}$, we have $o \notin \relint \conv (\bar{V} \cap H)$.
		\item[(iv)] A polytope $P$ is a pyramid if and only if at least one point of $\bar{V}$ coincides with
		the origin $o \in \R^{n-d-1}$.
	\end{itemize}
\end{theorem}

We note that (ii) can be stated in a more general form: $F$ is a face of $P$ if, and only if, for the corresponding co-face $\bar{F}$ of $P$, we have $o \in \inter \conv \bar{Z}$.

Before stating the result on $n=d+3$ vertices, recall that a $d$-polytope with $n$ vertices is \emph{cyclic}, if it is combinatorially equivalent to the convex hull of $n$ points on the moment curve $\gamma(t) = (t,t^2,\ldots, t^d)$, $t \in \R$.

\begin{theorem}[\cite{gho-langi 2}]\label{thm:dplus3}
	Let $P \in \P_d(d+3)$ satisfy Property Z. If $P$ is even, assume that $P$ is not cyclic. Then
	$P=\conv \{ P_1 \cup P_2 \cup P_3\}$, where $P_1$, $P_2$ and $P_3$ are regular simplices inscribed in $\Sph^{d-1}$ and contained in three mutually orthogonal linear subspaces of $\R^d$.
	Furthermore:
	\begin{itemize}
		\item If $d$ is odd and $P$ has maximal volume over $\P_d(d+3)$, then the dimensions of $P_1$, $P_2$ and $P_3$ are $\lfloor d/3 \rfloor$ or
		$\lceil d/3 \rceil$. In particular, in this case we have
		\[
		\left( v_d(d+3) = \right) \vol_d(P) = \frac{1}{d!} \cdot \prod_{i=1}^3 \frac{(k_i+1)^{\frac{k_i+1}{2}}}{k_i^\frac{k_i}{2}},
		\]
		where $k_1+k_2+k_3 = d$ and for every $i$, we have $k_i \in \left\{ \lfloor \frac{d}{3} \rfloor,  \lceil \frac{d}{3} \rceil \right\}$.
		\item The same holds if $d$ is even and $P$ has maximal volume over the family of \emph{not cyclic} elements of $\P_d(d+3)$.
	\end{itemize}
\end{theorem}

\begin{remark}\label{rem:cyclperm}
	Let $d=2m$ be even and $P \in \P_d(d+3)$ be a cyclic polytope satisfying Property Z.
	Then we need to examine the case that $\bar{P}$ is the vertex set of a regular $(2m+3)$-gon.
	Let the vertices of $\bar{P}$ be $\bar{p}_i$, $i=1,2,\ldots,2m+3$ in counterclockwise order.
	Applying the method of the proof of Theorem~\ref{thm:dplus3}, one can deduce that for every $i$, we have that $|p_{i-m-1}-p_i| = |p_{i+m+1}-p_i|$.
	On the other hand, for any other pair of vertices the conditions of Lemma~\ref{lem:note2} are not satisfied.
\end{remark}

In the light of Theorem~\ref{thm:dplus3}, it seems interesting to find the maximum volume cyclic polytopes in $\P_d(d+3)$, with $d$ even.
With regard to Remark~\ref{rem:cyclperm}, it is not unreasonable to consider the possibility that the answer for this question is a polytope $P =\conv \{p_i : i=1,2,\ldots,d+3 \}$ having a certain cyclic symmetry (if at all it is possible), namely that for any integer $k$, the value of $|p_{i+k}-p_i|$ is independent from $i$.

The following observation can be found both in \cite{KaibelWassmer}, or, as an exercise, in \cite{ziegler}.

\begin{remark}
	Let $d \geq 2$ be even, and $n \geq d+3$.
	Let
	\[
	C_d(n) = \sqrt{\frac{2}{d}} \conv \left\{ \left( \cos \frac{i\pi}{n}, \sin \frac{i\pi}{n},  \cos \frac{2i\pi}{n}, \ldots, \cos \frac{di\pi}{2n}, \sin \frac{di\pi}{2n} \right):i=0,1,\ldots,n-1 \right\}.
	\]
	Then $C(n,d)$ is a cyclic $d$-polytope inscribed in $\Sph^{d-1}$, and $\Sym(C_d(n)) = D_n$.
\end{remark}

It can be shown that for $d=4,6$ the only ``symmetric'' representations of a cyclic $d$-polytope with $d$ even and $n=d+3$ are those congruent to $C_d(d+3)$. Using the concepts of \emph{L\"owner ellipsoid} it can be proved the following theorem:

\begin{theorem}[\cite{gho-langi 2}]\label{thm:symmetric}
	Let $P \in \C_d(d+3)$ be a cyclic polytope, where $d=4$ or $d=6$, and let $V(P)=\{p_i:i=1,2,\ldots,d+3 \}$.
	If, for every value of $k$, $|p_{i+k}-p_i|$ is independent of the value of $i$, then $P$ is congruent to $C_d(d+3)$.
\end{theorem}

The above investigations raised a lot of questions and problems without answers. We collect some of them in the rest of this section.
\begin{problem}\label{prob:cyclic}
	Prove or disprove that in $\P_d(d+3)$, the cyclic polytopes with maximal volume are the congruent copies of $C_d(d+3)$.
	In particular, is it true for $C_4(7)$? Is it true that any cyclic polytope in $\P_d(d+3)$ satisfying Property Z is congruent to $C_d(d+3)$?
\end{problem}

A straightforward computation shows that $C_4(7)$ satisfies Property Z.
We note that in $\P_d(d+2)$, polytopes with maximal volume are cyclic, whereas in $\P_d(d+3)$, where $d$ is odd, they are not.
This leads to the following:

\begin{problem}\label{prob:whichone}
	Is it true that if $P \in \P_d(d+3)$, where $d$ is even, has volume $v_d(d+3)$, then $P$ is not cyclic?
\end{problem}

\begin{remark}\label{rem:notcyclic}
Let $P_4 \in \P_4(7)$ be the convex hull of a regular triangle and two diameters of $\Sph^3$, in mutually orthogonal linear subspaces.
Furthermore, let $P_6 \in \P_6(9)$ be the convex hull of three regular triangles, in mutually orthogonal linear subspaces.
One can check that
\[
\vol_4(P_4)= \frac{3}{4} = 0.43301\ldots > \vol_4(C_4(7)) = \frac{49}{192} \left( \cos \frac{\pi}{7} + \cos \frac{2\pi}{7} \right) = 0.38905\ldots.
\]
In addition,
\[
\vol_6(C_6(9))= \frac{7}{576} \sin \frac{\pi}{9} - \frac{7}{2880} \sin \frac{4\pi}{9} + \frac{7}{1152} \sin \frac{2\pi}{9} = 0.01697\ldots
\]
and
\[
\vol_6(P_6) = \frac{9\sqrt{3}}{640} = 0.02435\ldots > \vol_6(C_6(9)).
\]
This suggests that the answer for Problem~\ref{prob:whichone} is yes.
\end{remark}

\begin{remark}
	Using the idea of the proof of Theorem~\ref{thm:symmetric}, for any small value of $n$, it may identify the polytopes
	having $D_n$ as a subgroup of their symmetry groups. Nevertheless, it were unable to apply this method for general $n$, due to computational complexity. The authors \cite{gho-langi 2} carried out the computations for $5 \leq n \leq 9$, and obtained the following polytopes, up to homothety:
	\begin{itemize}
		\item regular $(n-1)$-dimensional simplex in $\R^{n-1}$ for every $n$,
		\item regular $n$-gon in $\R^2$ for every $n$,
		\item $C_4(n)$ with $n=6,7,8,9$ and $C_6(n)$ with $n=8,9$,
		\item regular cross-polytope in $\R^3$ and $\R^4$,
		\item the polytope $P_6$ in $\R^6$, defined in Remark~\ref{rem:notcyclic},
		\item the $3$-polytope $P$ with
		\[
		V(P)=\left\{ \left(1,0,0 \right),\left(-\frac{2}{3},-\frac{2}{3},\frac{1}{3} \right),  \left(0,1,0 \right), \left(\frac{1}{3},-\frac{2}{3},-\frac{2}{3} \right), \left(0,0,1 \right), \left(-\frac{2}{3},\frac{1}{3},-\frac{2}{3} \right)\right\} .
		\]
	\end{itemize}
\end{remark}

We note that, for $d$ odd, the symmetry group of a cyclic $d$-polytope with $n \geq d+3$ vertices is $\Z_2 \times \Z_2$ (cf. \cite{KaibelWassmer}).
Thus, the only cyclic polytopes in the above list are simplices and those homothetic to $C_d(n)$ for some values of $n$ and $d$.
This leads to the following question.

\begin{problem}
	Is it true that if, for some $n \geq d+3 \geq 5$, a cyclic polytope $P \in \P_d(n)$ satisfies $\Sym(P) = D_n$, then $P$ is congruent to $C_d(n)$?
\end{problem}

\section{Volume of the convex hull of two connecting bodies}

\subsection{On the volume function of the convex hull of two convex body}

Following the chronology, we have to start here with a result of F\'ary and R\'edei from 1950  (\cite{fary-redei}). They investigated the volume function defined on the convex hull of two convex bodies. He proved that if one of the bodies moves on a line with constant velocity then the volume of the convex hull is a convex function of the time (see Satz.4 in \cite{fary-redei}).
\begin{figure}[h]
\includegraphics[height=6cm]{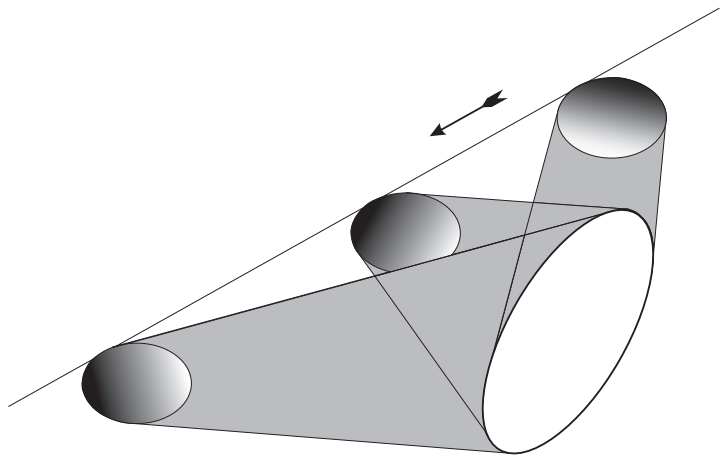}
\caption{The change of the convex hull.}
\label{fig:fary-redei}
\end{figure}
It was also proved later in \cite{rogers-shephard2}, and for convex polyhedra of dimension three in \cite{ahn}.

\begin{theorem}[\cite{fary-redei},\cite{rogers-shephard2},\cite{ahn}]\label{stm:convexity}
	The real valued function $g$ of the real variable $x$ defined by the fixed vector $t$ and the formula
	$$
	g(x):=\mathrm{ Vol }(\mathrm{ conv }(K\cup (K'+t(x))), \mbox{ where } t(x):=xt,
	$$
	is convex.
\end{theorem}

The nice proof in \cite{ahn} is based on the observation that the volume change function (by a translation in the direction of a line) can be calculated and it is an increasing function. Since it is also the derivative of $g$ we get that $g$ is convex. This calculation for the volume change can be done in the general case, too. Consider the shadow boundary of the convex hull $\mathrm{ conv }(K\cup (K'+t))$ with respect to the line of translation $t$. This is an $(d-2)$-dimensional topological manifold separating the boundary of
$\mathrm{ conv }(K\cup (K'+t))$ into two domains, the front and back sides of it, respectively. (The translation $t$ can be considered as a motion, hence the respective concepts of front and back sides can be regarded with respect to the direction of it.) Regarding a hyperplane $H$ orthogonal to $t$ the front side and  back side are graphs of functions over the orthogonal projection $X$ of $\mathrm{ conv }(K\cup (K'+t))$ onto $H$. Thus the volume change in $t$ can be calculated by the formula
$$
g'(t)=\lim\limits_{\varepsilon \rightarrow 0}\int_{X} (f^{t+\varepsilon}(X)-f^t(X))+\int_{X} (b^{t+\varepsilon}(X)-b^t(X)),
$$
where, at the moment $t$, $f^t$ and $b^t$ are the graphs of the front ad back sides, respectively. Since $X$ is independent from $t$ and for fixed $X$ the functions
$$
f^{t+\varepsilon}(x)-f^t(x) \mbox{ and } b^{t+\varepsilon}(x)-b^t(x)
$$
in $t$ are increasing and decreasing, respectively, we get that $g'$ is also increasing in $t$ implying that $g$ is convex.

As a corollary we get the following:
\begin{corollary}[see in \cite{gho 2}]
If we have two convex, compact bodies $K$ and $K'$ of the Euclidean space of dimension $n$ and they are moving uniformly on two given straight lines then the volume of their convex hull is a convex function of the time.
\end{corollary}	

\begin{remark} We emphasize that the statement of Theorem \ref{stm:convexity} is not true in hyperbolic space: Let $K$ be a segment and $K'$ be a point which goes on a line in the pencil of the rays ultraparalel to the line of the segment. Since the area function of the triangle defined by the least convex hull of $K$ and $K'$ is bounded (from below and also from above) it cannot be a convex function.
\end{remark}

Note that if the "bodies" are points the statement simplified to a proposition of absolute geometry which implies e.g. the existence of the normal transversal of two skew lines in the hyperbolic space.

There are several applications of Theorem \ref{stm:convexity}. In the paper of Hee-Kap Ahn, Peter Brass and Chan-Su Shin (see \cite{ahn}) the following result, based on Lemma \ref{lem:mainlemma} appears.
\begin{theorem}\label{thm.anh}(See Theorem 3 in \cite{ahn})
	Given two convex polyhedra $P$ and $Q$ in three-dimensional space, we can compute the translation vector $t$ of $Q$ that minimizes
	$\vol(\conv(P \cup (Q+ t)))$ in expected time
	$O(n^3 \log^4 n)$.
	The $d$-dimensional problem can be solved in expected time
	$O(n^{d+1-3/d} (log n)^{d+1})$.
\end{theorem}

In \cite{fary-redei}, F\'ary and R\'edey introduced the concepts of inner symmetricity (or outer symmetricity) of a convex body with the ratio (or inverse ratio) of the maximal (or minimal) volumes of the centrally symmetric bodies inscribed in (or circumscribed about) the given body. Using the mentioned Theorem \ref{stm:convexity} (and also its counterpart on the concavity of the volume function of the intersection of two bodies one of them from which moving on a line with constant velocity), they determined the inner symmetricity (and also the outer symmetricity) of a simplex (see Satz 5., resp. Satz 6. in \cite{fary-redei}). It has been proved that if $S$ is a simplex of dimension $n$ then its inner symmetricity $c_{\star}(S)$ is equal to
\begin{equation}
c_{\star}(S)=\frac{1}{(n+1)^n}\sum\limits_{0\leq \nu\leq \frac{n+1}{2}}(-1)^\nu\binom {n+1}{\nu}(n+1-2\nu)^n.
\end{equation}
On outer symmetricity $c^{\star}(S)$ of a simplex they proved that it is equal to
\begin{equation}
c^{\star}(S)=\frac{1}{\binom {n}{n_0}},
\end{equation}
where $n_0=n/2$ if $n$ is even and $n_0=(n-1)/2$ if $n$ is an odd number.
The above values attain when we consider the volume of the intersection (or the convex hull of the union) of $S$ with its centrally reflected copy $S_O$ (taking the reflection at the centroid $O$ of $S$).

\'A. G.Horv\'ath and Zs. L\'angi in \cite{gho-langi} introduced the following quantity.
\begin{defi}
	For two convex bodies $K$ and $L$ in $\R^d$, let
	\[
	c(K,L)=\max\left\{ \vol_d (\conv (K'\cup L')) : K' \cong K, L' \cong L \mbox{ and } K'\cap L'\neq\emptyset \right\}.
	\]
	Furthermore, if $\S$ is a set of isometries of $\R^d$, we set
	\[
	c(K|\S)= \frac{1}{\vol(K)} \max \left\{ \vol_d( \conv (K \cup K' )) : K \cap K' \neq \emptyset , K' = \sigma(K) \hbox{ for some } \sigma \in \S \right\} .
	\]
\end{defi}

A quantity similar to $c(K,L)$ was defined by Rogers and Shephard \cite{rogers-shephard2}, in which congruent copies were replaced by translates.
It has been shown that the minimum of $c(K|\S)$, taken over the family of convex bodies in $\R^d$, is its value for a $d$-dimensional Euclidean ball, if $\S$ is the set of translations or that of reflections about a point.
Nevertheless, their method, approaching a Euclidean ball by suitable Steiner symmetrizations and showing that during this process the examined quantities do not increase, does not characterize the convex bodies for which the minimum is attained;
they conjectured that, in both cases, the minimum is attained only for ellipsoids (cf. p. 94 of \cite{rogers-shephard2}).
We note that the method of Rogers and Shephard \cite{rogers-shephard2} was used also in \cite{macbeath}.
The results of the mentioned work based on the concept of \emph{linear parameter system of convex sets} and such a generalization of Theorem \ref{stm:convexity} which has interest on its own-right, too.

\begin{defi}\label{def:linparsys}(\cite{rogers-shephard2})
Let $I$ be an arbitrary index set, with each member $i$ of which is associated a point $a_i$ in $d$-dimensional space, and a real number $\lambda_i$, where the sets $\{a_i\}_{i\in I}$  and $\{\lambda_i\}_{i\in I}$ are each bounded. If $e$ is a fixed point and $t$ is any real number, $A(t)$ denotes the set of points
$$
\{a_i+t\lambda_ie\}_{i\in I},
$$
and $C(t)$ is the least convex cover of this set of points, then the system of convex sets $C(t)$ is called a \emph{linear parameter system}.
\end{defi}
The authors proved (see Lemma 1 in \cite{rogers-shephard2}) that the volume $V(t)$ of the set $C(t)$ of a linear parameter system is a convex function of $t$. They noted  that this result should be contrasted with that for a linear system of convex bodies as defined by Minkowski, where the $d$-th root of the volume of the body with parameter $t$ is a concave function of $t$ in its interval of definition.

In this paper we prove the following results:
\begin{theorem}\label{thm:rogers-shepard1}
Let $H$ and $K$ be two bodies and denote by $C(H,K)$ the least convex cover of the union of $H$ and $K$. Furthermore let $V^{\star}(H,K)$ denote the maximum, taken over all point $x$ for which the intersection $H\cap (K+x)$ is not empty, of the volume $\vol_d(C(H,K+x))$ of the set $C(H,K+x)$. Then $V^{\star}(H,K)\geq V^{\star}(SH,SK)$, where $SH$ denotes the closed $d$-dimensional sphere with centre at the origin and with volume equal to that of $H$.
\end{theorem}
\begin{theorem}\label{thm:rogers-shepard2}
If $K$ is a convex body in $d$-dimensional space, then
	$$
	1+\frac{2J_{d-1}}{J_d}\leq \frac{\vol_d(R^\star K)}{\vol_d(K)}\leq 2^d,
	$$
	where $J_d$ is the volume of the unit sphere in $d$-dimensional space, $R^\star K$ is the number to maximize with respect to a point $a$ of $K$ the volumes of the least centrally symmetric convex body with centre $a$ and containing $K$. Equality holds on the left, if $K$ is an ellipsoid; and on the right, if, and only if, $K$ is a simplex.
\end{theorem}
\begin{theorem}\label{thm:rogers-shepard3}
If $K$ is centrally symmetric body in $d$-dimensional space, then
	$$
	1+\frac{2J_{d-1}}{J_d}\leq \frac{\vol_d(R^\star K)}{\vol_d(K)}\leq 1+d
	$$
	Equality holds on the left if $K$ is an ellipsoid, and on the right if $K$ is any centrally symmetric double-pyramid on a convex base.
\end{theorem}
\begin{theorem}\label{thm:rogers-shepard4}
If $K$ is a convex body in $d$-dimensional space, then
		$$
		1+\frac{2J_{d-1}}{J_d}\leq \frac{\vol_d(T^\star K)}{\vol_d(K)}\leq 1+d,
		$$
	where $T^\star K$ denotes the so-called \emph{translation body} of $K$. This is the body for which the volume of $K\cap (K+x)\neq \empty$ and the volume of $C(K,K+x)$ is maximal one. Equality holds on the left if $K$ is an ellipsoid, and on the right if $K$ is a simplex.
\end{theorem}
\begin{theorem}\label{thm:rogers-shepard5}	
	Let $K$ be a convex body in $d$-dimensional space. Then
	there is a direction such that the volume of each cylinder $Z$, circumscribed
	to $K$, with its generators in the given direction, satisfies
	$$
	\frac{\vol_d(Z)}{\vol_d(K)}\geq \frac{2J_{d-1}}{J_d}.
	$$
\end{theorem}

It can be seen that these statements connect with the problem to determine the number $c(K|\S)$ defined in Definition 3.1. In fact, G. Horv\'ath and L\'angi (in \cite{gho-langi}) treated these problems in a more general setting. Let $c_i(K)$ be the value of $c(K|\S)$, where $\S$ is the set of reflections about the $i$-flats of $\R^d$, and $i=0,1,\ldots,d-1$.
Similarly, let $c^{tr}(K)$ and $c^{co}(K)$ be the value of $c(K|\S)$ if $\S$ is the set of translations and that of all the isometries, respectively.
In \cite{gho-langi} the authors examined the minima of these quantities. In particular, in Theorem~\ref{thm:translate}, was given another proof that the minimum of $c^{tr}(K)$, over the family of convex bodies in $\Re^n$, is its value for Euclidean balls, and it was shown also that the minimum is attained if, and only if, $K$ is an ellipsoid. This verifies the conjecture in \cite{rogers-shephard2} for translates.

Presented similar results about the minima of $c_1(K)$ and $c_{d-1}(K)$, respectively. In particular, the authors proved that, over the family of convex bodies, $c_1(K)$ is minimal for ellipsoids, and $c_{n-1}(K)$ is minimal for Euclidean balls. The first result proves the conjecture of Rogers and Shephard for copies reflected about a point.

During the investigation, $\K_d$ denotes the family of $d$-dimensional convex bodies. For any $K \in \K_d$ and $u \in \Sph^{n-1}$, $K | u^\perp$ denotes the orthogonal projection of $K$ into the hyperplane passing through the origin $o$ and perpendicular to $u$. The \emph{polar} of a convex body $K$ is denoted by $K^\circ$. The denotation $J_d$ of the paper \cite{rogers-shephard2} we are changing to the more convenient one $v_d$.)

The propositions are the followings:

\begin{theorem}\label{thm:translate}
	For any $K\in \K_d$ with $d\geq 2$, we have $c^{tr}(K) \geq 1 + \frac{2v_{d-1}}{v_d}$ with equality if, and only if, $K$ is an ellipsoid.
\end{theorem}

We remark that a theorem related to Theorem~\ref{thm:translate} can be found in \cite{MM06}. More specifically, Theorem 11 of \cite{MM06} states that
for any convex body $K \in \K_d$, there is a direction $u \in \Sph^{d-1}$ such that, using the notations of Theorem~\ref{thm:translate},
$d_K(u) \vol_{d-1}(K|u^\perp) \geq \frac{2v_{d-1}}{v_d}$, and if for any direction $u$ the two sides are equal, then $K$ is an ellipsoid.

If, for a convex body $K \in \K_d$, we have that $\vol_d (\conv ((v+K) \cup (w+K)))$ has the same value for any touching pair of translates,
let us say that $K$ satisfies the \emph{translative constant volume property}.
The characterization of the plane convex bodies with this property can be found also in this paper.
Before formulating the result, we recall that a $2$-dimensional $o$-symmetric convex curve is a Radon curve, if, for the convex hull $K$ of
a suitable affine image of the curve, it holds that $K^\circ$ is a rotated copy of $K$ by $\frac{\pi}{2}$ (cf. \cite{MS06}).

\begin{theorem}\label{thm:pointreflection}
	For any plane convex body $K \in \K_2$ the following are equivalent.
	\begin{itemize}
		\item[(1)] $K$ satisfies the translative constant volume property.
		\item[(2)] The boundary of $\frac{1}{2}(K-K)$ is a Radon curve.
		\item[(3)] $K$ is a body of constant width in a Radon norm.
	\end{itemize}
\end{theorem}

In two situations we have more precise results, respectively. The first case is when the examined body is a centrally symmetric one, and the other one when it is symmetric with respect to a hyperplane. The authors proved the following two theorems:

\begin{theorem}\label{thm:centralsymmetry}
	For any $K \in \K_d$ with $d \geq 2$, $c_1(K) \geq 1 + \frac{2v_{d-1}}{v_d}$, with equality if, and only if, $K$ is an ellipsoid.
\end{theorem}

\begin{theorem}\label{thm:hyperplane}
	For any $K \in \K_d$ with $d \geq 2$, $c_{d-1}(K) \geq 1+\frac{2v_{d-1}}{v_d}$, with equality if, and only if, $K$ is a Euclidean ball.
\end{theorem}

Finally, let $\P_m$ denote the family of convex $m$-gons in the plane $\R^2$.
It is a natural question to ask about the minima of the quantities defined in the introduction over $\P_m$.
More specifically, we set
\begin{eqnarray*}
	t_m & = & \min \{ c^{tr}(P) : P \in \P_m \};\\
	p_m & = & \min \{ c_0(P) : P \in \P_m\};\\
	l_m & = & \min \{ c_1(P) : P \in \P_m \}.
\end{eqnarray*}
On these numbers the following results were shown:
\begin{theorem}\label{thm:discrete}
		\begin{itemize}
		\item[(1)] $t_3 = t_4 = 3$ and $t_5 = \frac{25+\sqrt{5}}{10}$. Furthermore, $c^{tr}(P) = 3$ holds for any triangle and quadrangle, and if $c^{tr}(P) = t_5$ for some $P \in \P_5$, then $P$ is affine regular pentagon.
		\item[(2)] $p_3 = 4$, $p_4 = 3$ and $p_5 = 2 + \frac{4 sin \frac{\pi}{5}}{5}$. Furthermore, in each case, the minimum is attained only for affine regular polygons.
		\item[(3)] $l_3 = 4$ and $l_4 =3$.  Furthermore, among triangles, the minimum is attained only for regular ones, and among quadrangles for rhombi.
	\end{itemize}
\end{theorem}

\begin{conjecture}
	Let $d \geq 2$ and $0 < i < d-1$. Prove that, for any $K \in \K_d$, $c_i(K) \geq 1 + \frac{2v_{d-1}}{v_d}$. Is it true that equality holds only for Euclidean balls?
\end{conjecture}

The maximal values of $c^{tr}(K)$ and $c_0(K)$, for $K \in \K_d$, and the convex bodies for which these values are attained, are determined in \cite{rogers-shephard2}. Using a suitable simplex as $K$, it is easy to see that the set $\{c_i(K) : K \in \K_n \}$ is not bounded from above for
$i = 1,\ldots,n-1$. This readily yields the same statement for $c^{co}(K)$ as well.
On the other hand, from Theorem~\ref{thm:hyperplane} we obtain the following.

\begin{remark}
	For any $K \in \K_n$ with $n \geq 2$, we have $c^{co}(K) \geq 1+\frac{2v_{n-1}}{v_n}$, with equality if, and only if, $K$ is a Euclidean ball.
\end{remark}

In Theorem~\ref{thm:pointreflection} it was proved that in the plane, a convex body satisfies the translative equal volume property if, and only if, it is of constant width in a Radon plane. It is known (cf. \cite{AB89} or \cite{MS06}) that for $d \geq 3$, if every planar section of a normed space is Radon, then the space is Euclidean; that is, its unit ball is an ellipsoid. This motivates the conjecture:

\begin{conjecture}
	Let $d \geq 3$. If some $K \in \K_d$ satisfies the translative equal volume property, then $K$ is a convex body of constant width in a Euclidean space.
\end{conjecture}

Furthermore, we remark that the proof of Theorem~\ref{thm:pointreflection} can be extended, using the Blaschke-Santal\'o inequality, to prove Theorems~\ref{thm:translate} and \ref{thm:centralsymmetry} in the plane.
Similarly, Theorem~\ref{thm:hyperplane} can be proved by a modification of the proof of Theorem~\ref{thm:translate}, in which we estimate the volume of the polar body using the width function of the original one, and apply the Blaschke-Santal\'o inequality.

Like in \cite{rogers-shephard2}, Theorems~\ref{thm:translate} and \ref{thm:hyperplane} yield information about circumscribed cylinders.
Note that the second corollary is a strengthened version of Theorem 5 in \cite{rogers-shephard2}.

\begin{corollary}\label{cor:rightcylinder}
	For any convex body $K \in \K_d$, there is a direction $u \in \Sph^{d-1}$ such that the right cylinder $H_K(u)$, circumscribed about $K$
	and with generators parallel to $u$ has volume
	\begin{equation}\label{eq:rightcylinder}
	\vol(H_K(u)) \geq \left( 1 + \frac{2v_{d-1}}{v_d} \right) \vol_d(K).
	\end{equation}
	Furthermore, if $K$ is not a Euclidean ball, then the inequality sign in (\ref{eq:rightcylinder}) is a strict inequality.
\end{corollary}

\begin{corollary}\label{cor:eachcylinder}
	For any convex body $K \in \K_d$, there is a direction $u \in \Sph^{d-1}$ such that any cylinder $H_K(u)$, circumscribed about $K$
	and with generators parallel to $u$, has volume
	\begin{equation}\label{eq:eachcylinder}
	\vol(H_K(u)) \geq \left( 1 + \frac{2v_{d-1}}{v_d} \right) \vol_d(K).
	\end{equation}
	Furthermore, if $K$ is not an ellipsoid, then the inequality sign in (\ref{eq:eachcylinder}) is a strict inequality.
\end{corollary}

Let $P_m$ be a regular $m$-gon in $\R^2$.

\begin{problem}
	Prove or disprove that for any $m \geq 3$,
	\[
	t_m = c^{tr}(P_m), \quad p_m = c_0(P_m), \quad \hbox{and}  \quad l_m = c_1(P_m).
	\]
	Is it true that for $t_m$ and $p_m$, equality is attained only for affine regular $m$-gons, and for $l_m$, where $m \neq 4$, only for regular $m$-gons?
\end{problem}

\subsection{Simplices in the $3$-space}

\'A. G.Horv\'ath in \cite{gho} examined $c(K,K)$ in the special case that $K$ is a regular tetrahedron and the two congruent copies have the same centre. It has been proved the following theorem.
\begin{theorem}
	The volume of the convex hull of two congruent regular triangles
	with a common center is maximal if and only if their planes are
	orthogonal to each other and one of their vertices are opposite
	position with respect to the common center $O$.
\end{theorem}

The proof is based on some exact formulas, which can be extended to the nonregular case of triangles, too.

On regular tetrahedra was proved a theorem in that case when all of the spherical triangles contain exactly one from the vertices of the other tetrahedron and changing the role of the tetrahedra we also get it (so when the two tetrahedra are in \emph{ dual position}). We remark that in a dual position the corresponding spherical edges of the two tetrahedra are crossing to each other, respectively. In this case it has been proved that

\begin{theorem}
	The value $v=\frac{8}{3\sqrt{3}}r^3$ is an upper bound for the volume of the convex hull of two regular tetrahedra are in dual position. It is attained if and only if the eight vertices of the two tetrahedra are the vertices of a cube inscribed in the common circumscribed sphere.
\end{theorem}

This paper considered the proof of that combinatorial case when two domains contain two vertices, respectively. The following statement were proved:

\begin{statement}
	Assume that the closed regular spherical simplices $S(1,2,3)$ and $S(4,2,3)$ contains the vertices $2',4'$ and $1',3'$, respectively. Then the two tetrahedra are the same.
\end{statement}

In the paper \cite{gho 3} the author closed this problem using a generalization of the icosahedron inequality of L. Fejes-T\'oth. It has been shown the general statement:

\begin{theorem}
	Consider two regular tetrahedra inscribed in the unit sphere. The maximal volume of the convex hull $P$ of the eight vertices is the volume of the cube $C$ inscribed in the unit sphere, so
	$$
	\vol_3(P)\leq \vol_3(C)=\frac{8}{3\sqrt{3}}.
	$$
\end{theorem}

The paper \cite{gho 2} investigates also connecting simplices. It is assumed that the used set of isometries $\S$ consists only reflections at such hyperplanes $H$ which intersect the given simplex $S$. (Hence the convex hull function considered only on the pairs of the simplex and its reflected copy at a hyperplane intersecting it.) Explicitly wrote the relative volume of the convex hull of the simplices and gave upper bounds on it. The number $c(S,S^H)$ for the regular simplex determined explicitly. The following lemma plays a fundamental role in the investigations.

\begin{lemma}\label{lem:mainlemma}
	If $K$ and $K'$ give a maximal value for $c_{K,K'}$
	then the intersection $K\cap K'$ is an extremal point of each of the bodies.
\end{lemma}

To formulate the results we introduce some new notation. Assume that the intersecting simplices $S$ and $S_H$ are reflected copies of each other at the hyperplane $H$. Then $H$ intersects each of them in the same set. By the Lemma \ref{lem:mainlemma} we have that the intersection of the simplices in an optimal case is a common vertex. Let $s_0\in H$ and $s_i\in H^{+}$ for $i\geq 1$. We imagine that $H$ is horizontal and $H^+$ is the upper half-space. Define the \emph{upper side of $S$} as the union of those facets in which a ray orthogonal to $H$ and terminated in a far point of $H^+$ is first intersecting with $S$.  The volume of the convex hull is the union of those prisms which are based on the orthogonal projection of a facet of the simplex of the upper side. Let denote $F_{i_1},\cdots, F_{i_k}$ the facet-simplex of the upper side, $F'_{i_1},\cdots, F'_{i_k}$ its orthogonal projections on $H$ and $u_{i_1}, \cdots, u_{i_k}$ its respective unit normals, directed outwardly. We also introduce the notation $s=\sum\limits_{i=0}^{d}s_i=\sum\limits_{i=1}^{d}s_i$. Now we have

\begin{statement}
	$$
	\frac{1}{\mathrm{ vol }_{d}(S)}\mathrm{ vol }_d(\mathrm{ conv }(S,S^H))=2d\sum\limits_{l=1}^k\frac{\langle u_{i_l}, u\rangle \langle u, s-s_{i_l}\rangle}{|\langle u_{i_l}, (d+1)s_{i_l}-s\rangle |}.
	$$
\end{statement}

It can be solved the original problem in the case of the regular simplex. Denote the Euclidean norm of a vector $x$ by $\|x\|$.

\begin{theorem}
	If $S$ is  the regular simplex of dimension $n$, then
	$$
	c(S,S^H):=\frac{1}{\mathrm{ vol }_{d}(S)}\mathrm{ vol }_d(\mathrm{ conv }(S,S^H))=2d,
	$$
	attained only in the case when $u=u_0=\frac{s}{\|s\|}$.
\end{theorem}

We note that the result of the case of reflection at a hyperplane gives an intermediate value between the results corresponding to translates and point reflections. The part of the previous proof corresponding to the case of a single upper facet can be extended to a general simplex, too. Let $G$ denote the Gram matrix of the vector system $\{s_1,\ldots,s_n\}$, defined by the product $M^TM$, where $M=[s_1,\cdots,s_n]$ is the matrix with columns $s_i$. In the following theorem we use the notation $\|\cdot\|_1$ for the $l_1$ norm of a vector or a matrix, respectively.

\begin{theorem}
	If the only upper facet is $F_0$ with unit normal vector $u_0$, then we have the inequality
	$$
	\frac{1}{\mathrm{ vol }_{d}(S)}\mathrm{ vol }_d(\mathrm{ conv }(S,S^H))\leq d\left(1+\frac{\|s\|}{\langle u_0,s\rangle }\right)=
	$$
	$$
	=\left(d+\sqrt{\left\|(1,\ldots,1)G^{-1}\right\|_1}\left\|M(1,\ldots,1)\right\|\right).
	$$
	Equality is attained if and only if the normal vector $u$ of $H$ is equal to $\frac{u_0+s'}{\|u_0+s'\|}$, where $s'=\frac{s}{\|s\|}$ is the unit vector of the direction of $s$.
\end{theorem}

We remark that for a regular simplex we get back the previous theorem, since
$$
G=\left(
\begin{array}{cccc}
1 & \frac{1}{2} & \cdots & \frac{1}{2} \\
\frac{1}{2} & 1 & \cdots & \frac{1}{2} \\
\vdots & \vdots & \vdots & \vdots \\
\frac{1}{2} & \cdots & \frac{1}{2} & 1 \\
\end{array}
\right) \mbox{ and }
G^{-1}=\left(
\begin{array}{cccc}
\frac{2d}{d+1} & -\frac{2}{d+1} & \cdots & -\frac{2d}{d+1} \\
-\frac{2}{d+1} & \frac{2d}{d+1} & \cdots & -\frac{2}{d+1} \\
\vdots & \vdots & \vdots & \vdots \\
-\frac{2}{d+1} & \cdots & -\frac{2}{d+1} & \frac{2d}{d+1} \\
\end{array}
\right),
$$
implying that
$$
d+\sqrt{\left\|(1,\ldots,1)G^{-1}\right\|_1}\left\|M(1,\ldots,1)\right\|=d+\sqrt{\frac{2d}{d+1}}\sqrt{\frac{d(d+1)}{2}}=2d.
$$

\end{document}